\documentclass[12pt,oneside,reqno]{amsart}
\usepackage{txfonts}
\usepackage{bbm}
\usepackage{amsmath}
\usepackage{graphicx}
\usepackage{mathrsfs}
\usepackage{stmaryrd}
\usepackage{amsfonts}
\usepackage{enumerate,amsmath,amssymb,amsthm}

\numberwithin{equation}{section}

\newcommand{\be}{\begin{eqnarray}}
\newcommand{\ee}{\end{eqnarray}}
\newcommand{\ce}{\begin{eqnarray*}}
\newcommand{\de}{\end{eqnarray*}}
\newtheorem{theorem}{Theorem}[section]
\newtheorem{lemma}[theorem]{Lemma}
\newtheorem{remark}[theorem]{Remark}
\newtheorem{definition}[theorem]{Definition}
\newtheorem{proposition}[theorem]{Proposition}
\newtheorem{Examples}[theorem]{Example}
\newtheorem{corollary}[theorem]{Corollary}

\def\eps{\varepsilon}

\def\p{\partial}

\def\[{{\Big[}}
\def\]{{\Big]}}
\def\<{{\langle}}
\def\>{{\rangle}}
\def\({{\Big(}}
\def\){{\Big)}}

\def\bx{{\mathbf{x}}}

\def\dif{{\mathord{{\rm d}}}}

\def\min{{\mathord{{\rm min}}}}

\def\no{\nonumber}
\def\={&\!\!=\!\!&}
\def\bt{\begin{theorem}}
\def\et{\end{theorem}}
\def\bl{\begin{lemma}}
\def\el{\end{lemma}}
\def\br{\begin{remark}}
\def\er{\end{remark}}
\def\bd{\begin{definition}}
\def\ed{\end{definition}}
\def\bp{\begin{proposition}}
\def\ep{\end{proposition}}
\def\bc{\begin{corollary}}
\def\ec{\end{corollary}}
\def\bx{\begin{Examples}}
\def\ex{\end{Examples}}

\def\cL{{\mathcal L}}

\def\cS{{\mathcal S}}

\def\mD{{\mathbb D}}
\def\mE{{\mathbb E}}

\def\mH{{\mathbb H}}

\def\mN{{\mathbb N}}

\def\mQ{{\mathbb Q}}
\def\mR{{\mathbb R}}
\def\mS{{\mathbb S}}

\def\mU{{\mathbb U}}

\def\mW{{\mathbb W}}

\def\sB{{\mathscr B}}
\def\sC{{\mathscr C}}

\def\sF{{\mathscr F}}

\def\geq{\geqslant}
\def\leq{\leqslant}

\def\div{\mathord{{\rm div}}}

\def\v{{\mathrm v}}
\def\e{{\mathrm e}}

\allowdisplaybreaks

\begin{document}

\title{Fundamental solution of kinetic Fokker-Planck operator with anisotropic nonlocal dissipativity}

\date{}
\author{Xicheng Zhang}

\address{
School of Mathematics and Statistics,
Wuhan University, Wuhan, Hubei 430072, P.R.China\\
Email: XichengZhang@gmail.com
 }

\begin{abstract}
By using the probability approach (the Malliavin calculus), we prove the existence of smooth fundamental solutions for degenerate kinetic Fokker-Planck equation with  anisotropic nonlocal dissipativity, where the dissipative term is the generator of an anisotropic L\'evy process and
the drift term is allowed to be cubic growth.
\end{abstract}

\maketitle \rm

\section{Introduction and Main Result}
Consider the following second order stochastic differential equation (SDE) in $\mR^d$:
\begin{align}
\frac{\dif^2 X_t}{\dif t^2}=-\nabla V(X_t)+\frac{\dif W_t}{\dif t}-\frac{\dif X_t}{\dif t},\ \ X_0=x,\label{SDE}
\end{align}
where $V(x):\mR^d\to\mR_+$ is a smooth function, and $(W_t)_{t\geq 0}$ is a $d$-dimensional standard Brownian motion.
In phase space $\mR^d_x\times\mR^d_{\v}$, the position and velocity vector field $(X_t,\dot X_t)$ solves the following degenerate SDE:
\begin{align}
\left\{\label{SDE00}
\begin{aligned}
\dif X_t&=\dot X_t\dif t,& X_0=x,\\
\dif\dot X_t&=-\nabla V(X_t)\dif t-\dot X_t\dif t+\dif W_t,& \dot X_0=\v.
\end{aligned}
\right.
\end{align}
The celebrated H\"ormander's hypoellipticity theorem tells us that $(X_t,\dot X_t)$ admits a smooth density $\rho_{x,\mathrm{v}}(t,x',\v')$ (cf. \cite{Ho, Ku-St2, No, Ma, Nu}).
Moreover, by It\^o's formula, one knows that $\rho_{t,\v}(t,x',\v')$ solves the following kinetic Fokker-Planck equation:
$$
\p_t\rho-\v'\cdot\nabla_{x'}\rho+\nabla V\cdot\nabla_{\v'}\rho=\Delta_{\v'}\rho-\div(\v'\rho).
$$
It is easy to check that the equilibrium of this equation is given by
$$
\rho_\infty(x,\v):=\exp\{-H(x,\v)\},\quad \mbox{where }\ \ \ H(x,\v):=\tfrac{|\v|^2}{2}+V(x).
$$
The rate of convergence to the equilibrium for the kinetic Fokker-Planck equation has been deeply studied in \cite{De-Vi, He-Ni, Vi, Gu-Wa}, etc.
Moreover, the stochastic flow property of SDE (\ref{SDE00}) was proven in \cite{Ar, Zh0}.

In this work, we shall consider equation (\ref{SDE}) with Brownian motion $(W_t)_{t\geq 0}$ replaced by
a L\'evy process $(L_t)_{t\geq 0}$ (for example, the cylindrical $\alpha$-stable process). 
More generally, we consider the following stochastic Hamiltonian system driven by L\'evy process:
\begin{align}
\left\{
\begin{aligned}
\dif X_t&=b_1(X_t,\dot X_t)\dif t,& X_0=x,\\
\dif\dot X_t&=b_2(X_t,\dot X_t)\dif t+\dif L_t,& \dot X_0=\v,
\end{aligned}
\right.\label{SDE0}
\end{align}
where $b=(b_1,b_2)$ is a smooth vector field on phase space $\mR^{d}_x\times\mR^d_\v$. The background about stochastic Hamiltonian system and
related Fokker-Planck equation is refereed to \cite{So}. From the microscopic viewpoint, stochastic equation (\ref{SDE0}) can be considered that
the motion of particles is perturbed by a ``discontinuous'' stochastic force.  We want to study the regularizing effect of L\'evy noise to the system.
It is well known that there are a lot of works devoting to the study of smooth densities for SDEs with jumps (see \cite{Bi,Bi-Ja-Gr,Pi,Ta,Is-Ku,Ca,Ku,Ba-Cl}, etc.).
Nevertheless, most of these works required that the jump noise is {\it non-degenerate}, and the main arguments are based upon developing an analogue of 
the Malliavin calculus for jump diffusions.

The main goal of the present paper is to prove that under some assumptions on $b$ and $L_t$, the solution $(X_t,\dot X_t)$ of SDE (\ref{SDE0})
still has a smooth density. When $b$ has bounded derivatives of all orders, $(\nabla_\v b_1)(\nabla_\v b_1)^*$ is uniform positive with
respect to $(x,\v)$, and $L_t$ is an isotropic $\alpha$-stable process, the smoothness of $\rho$ was proved in \cite{Zh2}. However, 
in real model such as stochastic oscillators, the nonlinear term $b$ is usually non-Lipschitz, and the L\'evy noise may be {\it anisotropic} 
as that each component of $L_t$ is independent.

Below, we first describe the noise $L_t$ following \cite{Ku}. Let $(L_t)_{t\geq 0}$ be a $d$-dimensional L\'evy process 
with the following form (called subordinated Brownian motion):
$$
L_t:=W_{S_t}=\Big(W^1_{S^1_t},\cdots,W^d_{S^d_t}\Big),
$$
where $S_t=(S^1_t,\cdots,S^d_t)$ is an independent $d$-dimensional $\mR^d_+$-valued L\'evy process with characteristic triple $(\vartheta,0,\nu_S)$ , more precisely, its Laplace transform is given by
\begin{align}
\mE(\e^{-z\cdot S_t})=\exp\left\{-t\vartheta\cdot z+\int_{\mR^d_+}(\e^{-z\cdot u}-1)\nu_S(\dif u)\right\},\label{EW111}
\end{align}
where $\vartheta\in\mR^d_+$ and the L\'evy measure $\nu_S$ satisfies
$$
\int_{\mR^d_+}(1\wedge |u|)\nu_S(\dif u)<\infty.
$$
In particular, each component $S^i_t$ is a subordinator (cf. \cite{Be,Sa}). By easy calculations, one can see that the characteristic function of $L_t$ is given by
\begin{align}
\mE \mathrm{e}^{\mathrm{i} z\cdot L_t}=\exp\left\{-t\sum_k\vartheta_k|z_k|^2+t\int_{\mR^d}(\mathrm{e}^{\mathrm{i}z\cdot y}
-1-\mathrm{i}z\cdot y1_{|y|\leq 1})\nu_L(\dif y)\right\},\label{Ch}
\end{align}
where $\nu_L$ is the L\'evy measure given by
\begin{align}
\nu_L(\Gamma)=\int_{\mR^d_+}\left(\int_\Gamma\frac{(2\pi)^{-d/2}}{(u_1\cdots u_d)^{\frac{1}{2}}} 
\mathrm{e}^{-(\frac{y_1^2}{2u_1}+\cdots+\frac{y_d^2}{2u_d})}\dif y_1\cdots\dif y_d\right)\nu_S(\dif u_1,\cdots,\dif u_d).\label{EW2}
\end{align}
Here we use the convention that if $u_i=0$ for some $i$, then the inner integral is calculated with respect to the degenerate Gaussian distribution.
In particular, $\nu_L$ may not be absolutely continuous with respect to the Lebesgue measure. Obviously, $\nu_L$ is a symmetric measure. 

Now we state the main assumptions on $(\vartheta,\nu_S)$ and  $b$:
\begin{enumerate}[{\bf(H$^{1}_{\nu_S}$)}]
\item Let $\phi:\mR_+\to\mR_+$ be defined by
\begin{align}
\phi(\eps):=\min_{i=1,\cdots, d}\left(\vartheta_i+\frac{1}{\e}\int_{|u|\leq\eps}u_i\nu_S(\dif u)\right).\label{Con2}
\end{align}
We assume that for some $\theta\in(0,1]$,
\begin{align}
\lim_{\eps\downarrow 0}\eps^{\theta-1}\phi(\eps)>0.\label{Con1}
\end{align}
\end{enumerate}
\begin{enumerate}[{\bf(H$^{2}_{\nu_S}$)}]
\item We assume that for any $p>0$,
\begin{align}
\int_{|u|>1}\e^{p|u|}\nu_S(\dif u)<\infty,\label{Con3}
\end{align}
which, by \cite[p.159, Theorem 25.3]{Sa},  is equivalent to
\begin{align}
\mE \e^{pS_t}<\infty.\label{ER1}
\end{align}
\end{enumerate}
\begin{enumerate}[{\bf(H$_b$)}]
\item Assume that there exists a Lyapunov function $H:\mR^d_x\times\mR^d_\v\to\mR_+$ with
\begin{align}
|\nabla_\v H|^2\leq C_1H,\ \ |\nabla^2_\v H|\leq C_2,\label{G1}
\end{align}
and such that for any $m\in\{0\}\cup\mN$ and some $q_m\geq 0$,
\begin{align}
b\cdot\nabla H\leq C_3H,\ |\nabla^m b|\leq C_m(H^{q_m}+1),\label{G2}
\end{align}
where $q_1\in[0,\frac{1}{2}]$. Moreover, 
\begin{align}
|\nabla_\v b|+|\nabla^2_\v b|+|\nabla^3_\v b|\leq C_4,\label{G3}
\end{align}
and for any row vector $u\in\mR^d$,
\begin{align}
|u\nabla_\v b_1(x,\v)|^2\geq C_5|u|^2.\label{G5}
\end{align}
\end{enumerate}
\br
Let $(S^i_t)_{i=1,\cdots,d}$ be independent $\alpha_i$-stable subordinators, where $\alpha_i\in(0,1)$. It is easy to check that (\ref{Con1}) holds with 
$\theta=\min(\alpha_1,\cdots,\alpha_d)$.
\er
\br
In the case of equation (\ref{SDE}), we can take 
$$
H(x,\v)=\frac{1}{2}|\v|^2+V(x),
$$
where $V\in C^\infty(\mR^d;\mR_+)$ satisfies that for any $m\in\mN$ and some $q_m\geq 0$,
$$
|\nabla^m V(x)|\leq C_m(V(x)^{q_m}+1),
$$
with $q_2=\frac{1}{2}$. In particular, $V(x)=|x|^4$ satisfies this assumption. Since for any $p>1$, compared with (\ref{ER1}),
it holds in general that (cf. \cite[p.168, Theorem 26.1]{Sa})
$$
\mE \e^{S_t^p}=\infty,
$$
we have to require $q_1\in[0,\frac{1}{2}]$ in (\ref{G2}) (see Theorem \ref{Th1} below).
\er

The main result of this paper is:
\bt\label{Main}
Under  {\bf(H$^{1}_{\nu_S}$)}, {\bf(H$^{2}_{\nu_S}$)} and {\bf(H$_b$)},
there exists a smooth density $\rho_{x,\v}(t,x',\v')$ to SDE (\ref{SDE0}) with bounded derivatives of all orders with respect to $x',\v'$ and such that
$$
\p_t\rho=\cL_{\v'}\rho+\div_{\v'} (b \rho),\quad r>0,
$$
with $\rho_{x,\v}(0,x',\v')=\delta_{x,\v}(x',v')$, where
$$
\cL_\v f(\v)=\mathrm{P.V.}\int_{\mR^d}(f(\v+y)-f(\v))\nu_L(\dif y)+\frac{1}{2}\sum_{k}(\p^2_k f)(\v)\vartheta_k,
$$
where $\mathrm{P.V.}$ stands for the Cauchy principal value. Moreover, there exist constants $\beta_1,\beta_2,\beta_3>0$ only depending on $d,\theta$ and
a positive continuous function $(x,\v)\mapsto C_{x,\v}$ such that for all $(t,(x',\v'),(x,\v))\in(0,1]\times(\mR^{d}_x\times\mR^{d}_\v)^2$,
\begin{align}
\rho_{x,\v}(t,x',\v')\leq C_{x,\v} \left(t^{-\beta_1}\left(1\wedge\frac{t^{\beta_2}}{(|x-x'|+|\v-\v'|)^{\beta_3}}\right)\right).\label{EE1}
\end{align}
\et
\br
If $b\in C^\infty_b(\mR^d)$, then the above $C_{x,\v}$ can be constant. In this case, if one only requires the existence of smooth density, 
then assumption {\bf(H$^{2}_{\nu_S}$)} can be dropped by using the same argument as in \cite[Section 3.3]{Zh2}.
\er
In order to prove this theorem, by taking regular conditional expectations with respect to $S_\cdot$, 
we shall regard the solution of SDE (\ref{SDE0}) as a Wiener functional, and then use the classical Malliavin calculus to prove Theorem \ref{Main}. 
Such an idea was first used by L\'endre \cite{Le}, and then in \cite{Ku, Zh1}. We also mention that a derivative formula of Bismut type and  the Harnack inequality
for SDEs driven by $\alpha$-stable processes were derived in \cite{Zh1} and \cite{Wa-Wa} following the same idea. It is quite interesting to have
an analytic proof of Theorem \ref{Main}.  It should be noticed that the L\'evy measure $\nu_L$ could be very singular. This leads to that the symbol of operator $\cL_\v$
$$
\Phi(\xi):=\int_{\mR^d}\(1-\e^{\mathrm{i}\<\xi,y\>}-1_{|y|\leq1}\<\xi,y\>\)\nu_L(\dif y)+\frac{1}{2}\sum_k\vartheta_k \xi_k^2
$$
may not be $C^1$-continuous differentiable on $\mR^d\setminus\{0\}$. 
Thus, the classical pseudo-differential operator theory seems not applicable (cf. \cite{Ho}). Below, we list some open questions for further studies:
\begin{itemize}
\item Can we prove the same result for multiplicative L\'evy noise?
\item Is it possible to remove the assumptions $q_1\in[0,\frac{1}{2}]$ in (\ref{G2}) and {\bf(H$^{2}_{\nu_S}$)}?
\item Is there a stationary distribution for stochastic Hamiltonian system (\ref{SDE0})? If yes, how about the rate of convergence as $t\to\infty$?
\end{itemize} 

This work is organized as follows: In Section 2, we prepare some notations and lemmas for later use. In particular, a Norris' type lemma is proven. 
In Section 3,  we prove some exponential moment estimate about the SDE driven by $W_{S_t}$ and with polynomial growth coefficients.
In Section 4, we calculate the Malliavin covariance matrix for the solution of SDE as a Wiener functional. In Section 5, 
we prove the smoothness of distributional density of a degenerate SDE driven by $W_{S_t}$. Meanwhile, we conclude the proof of Theorem \ref{Main}.
Before concluding this introduction, we collect some notations or conventions for later use.
\begin{itemize}
\item Write $\mR^d_+=[0,\infty)^d$ and $\mN_0=\{0\}\cup\mN$.
\item The inner product in Euclidean space is denoted by $\<x,y\>$ or $x\cdot y$.
\item For a vector $x=(x_1,\cdots,x_d)$, we write $|x|:=\left(\sum_i|x_i|^2\right)^{1/2}\sim\sum_i|x_i|$.
\item $C^\infty_0(\mR^d)$: The space of all smooth functions with compact support.
\item $\cS(\mR^d)$: The Schwardz space of rapidly decreasing  smooth functions.
\item $C^\infty_b(\mR^d)$: The space of all smooth bounded functions with bounded derivatives of all orders.
\item $C^\infty_p(\mR^d)$: The space of all smooth functions, which together with the derivatives of all orders are at most polynomial growth.
\item The asterisk $*$ denotes the transpose of a matrix or a column vector, or the dual operator.
\item $\nabla$ denotes the gradient operator, and $D$ the Malliavin derivative operator.
\item $C$ with or without index will denotes an unimportant positive constant.
\end{itemize}

\section{Preliminaries}

We first introduce the canonical space of subordinated Brownian motion $W_{S_t}$. Let $(\mW,\mH,\mu_\mW)$ be the classical Wiener space, i.e., $\mW$ is the space of all continuous functions from $\mR_+$ to $\mR^d$ with vanishing values at starting point $0$, $\mH\subset\mW$ is the Cameron-Martin space consisting of all
absolutely continuous functions with square integrable derivatives, $\mu_\mW$ is the Wiener measure so that the coordinate process
$$
W_t(w):=w_t
$$
is a $d$-dimensional standard Brownian motion. Let $\mS$ be the space of all c\`adl\`ag functions from $\mR_+$ to $\mR^d_+$ with
$\ell_0=0$, where each component is increasing. Suppose that $\mS$
is endowed with the Skorohod metric and
the probability measure $\mu_\mS$ so that the coordinate process
$$
S_t(\ell):=\ell_t=(\ell^1_t,\cdots,\ell^d_t)
$$
is a $d$-dimensional L\'evy process with Laplace transform (\ref{EW111}).
Consider the following product probability space
$$
(\Omega,\sF,P):=\Big(\mW\times \mS, \sB(\mW)\times\sB(\mS), \mu_\mW\times\mu_{\mS}\Big),
$$
and define for $(w,\ell)\in\mW\times \mS$,
$$
L_t(w,\ell):=w_{\ell_t}:=\Big(w_1(\ell_1(t)),\cdots, w_d(\ell_d(t))\Big).
$$
Then $(L_t)_{t\geq 0}$ is a L\'evy process with characteristic function (\ref{Ch}). We use the following filtration:
$$
\sF_t:=\sigma\{W_{S_s},S_s: s\leq t\}.
$$
Clearly, for $t>s$, $W_{S_t}-W_{S_s}$ and $S_t-S_s$ are independent of $\sF_s$.
\subsection{An exponential estimate of $S_t$}
The following estimate of exponential type about $S_t$ will play an important role in the proof of Theorem \ref{Main}.
\bl\label{Le1}
Let $f_t:\mR_+\to\mR^d_+$ be a continuous $\sF_t$-adapted process. For any $R,\eps,\delta>0$, we have
$$
P\left\{\int^{t\wedge\tau_R}_0f_s\cdot\dif S_s\leq\eps;\int^{t\wedge\tau_R}_0|f_s|\dif s>\delta\right\}
\leq \mathrm{e}^{1-\phi(\eps/R)\delta/\eps},
$$
where $\tau_R:=\inf\{t\geq 0: |f_t|>R\}$ and $\phi$ is defined by (\ref{Con2}).
\el
\begin{proof}
For $\lambda>0$, set
$$
g^\lambda_s:=\int_{\mR^d_+}(1-\mathrm{e}^{-\lambda f_s\cdot u})\nu_S(\dif u)
$$
and
$$
M^\lambda_t:=-\lambda\int^t_0f_s\cdot\dif S_s+\lambda\int^t_0f_s\cdot\vartheta\dif s+\int^t_0g^\lambda_s\dif s.
$$
Let $\mu(t,\dif u)$ be the Poisson random measure associated with $S_t$, i.e.,
$$
\mu(t,\Gamma):=\sum_{s\leq t}1_\Gamma(S_s-S_{s-}),\ \ \Gamma\in\sB(\mR^d_+).
$$
Let $\tilde \mu(t,\dif u)$ be the compensated Poisson random measure of $\mu(t,\dif u)$, i.e.,
$$
\tilde \mu(t,\dif u)=\mu(t,\dif u)-t\nu_S(\dif u).
$$
Then, by L\'evy-It\^o's decomposition (cf. \cite{Sa}), we can write
\begin{align}
S_t=t\left(\vartheta+\int_{|u|\leq 1}u\nu_S(\dif u)\right)+\int_{|u|\leq 1}u \tilde\mu(t,\dif u)+\int_{|u|>1}u \mu(t,\dif u),\label{EUY}
\end{align}
and so,
\begin{align*}
&\int^t_0f_s\cdot\dif S_s=\int^t_0 f_s\cdot\left(\vartheta+\int_{|u|\leq 1}u\nu_S(\dif u)\right)\dif s\\
&\qquad+\int^t_0\!\!\!\int_{|u|\leq 1}f_s\cdot u\tilde \mu(\dif s,\dif u)+\int^t_0\!\!\!\int_{|u|>1}f_s\cdot u\mu(\dif s,\dif u).
\end{align*}
By It\^o's formula (cf. \cite{Ap}), we have
\begin{align}
\mathrm{e}^{M^\lambda_t}=1+\int^t_0\!\!\!\int_{\mR^d_+}\mathrm{e}^{M^\lambda_{s-}}
[\mathrm{e}^{-\lambda f_s\cdot u}-1]\tilde \mu(\dif s,\dif u).\label{EQ2}
\end{align}
Since for any $x>0$,
$$
1-\mathrm{e}^{-x}\leq 1\wedge x,
$$
we have
$$
g^\lambda_s\leq\int_{\mR^d_+}(1\wedge(\lambda f_s\cdot u))\nu_S(\dif u)
$$
and
$$
M^\lambda_{t\wedge\tau_R}\leq \lambda\int^{t\wedge\tau_R}_0f_s\cdot\vartheta\dif s+\int^{t\wedge\tau_R}_0g^\lambda_s\dif s
\leq tR|\vartheta|+t\int_{\mR^d_+}(1\wedge(\lambda R |u|))\nu_S(\dif u).
$$
Hence, by (\ref{EQ2}) we have
$$
\mE \mathrm{e}^{M^\lambda_{t\wedge\tau_R}}=1.
$$
On the other hand, since for any $\kappa\in(0,1)$ and $x\leq-\log k$,
$$
1-\mathrm{e}^{-x}\geq \kappa x,
$$
letting $\kappa=\frac{1}{\e}$, we have for $s\leq \tau_R$,
\begin{align*}
\lambda f_s\cdot\vartheta+g^\lambda_s&\geq\lambda f_s\cdot\vartheta+\int_{|u|\leq\frac{1}{\lambda R}}(1-\mathrm{e}^{-\lambda f_s\cdot u})\nu_S(\dif u)\\
&\geq\lambda f_s\cdot\vartheta+\frac{1}{\e}\int_{|u|\leq\frac{1}{\lambda R}}(\lambda f_s\cdot u)\nu_S(\dif u)\\
&=\lambda f_s\cdot\left(\vartheta+\frac{1}{\e}\int_{|u|\leq\frac{1}{\lambda R}}u\nu_S(\dif u)\right)\\
&\geq\lambda\phi(1/(\lambda R))|f_s|,
\end{align*}
where $\phi$ is defined by (\ref{Con2}). Thus,
\begin{align*}
&\left\{\int^{t\wedge\tau_R}_0f_s\cdot\dif S_s\leq\eps;\int^{t\wedge\tau_R}_0|f_s|\dif s>\delta\right\}\\
&\quad\subset\left\{\mathrm{e}^{M^\lambda_{t\wedge\tau_R}}\geq \mathrm{e}^{-\lambda\eps
+\int^{t\wedge\tau_R}_0(\lambda f_s\cdot\vartheta+g^\lambda_s)\dif s};
\int^{t\wedge\tau_R}_0(\lambda f_s\cdot\vartheta+g^\lambda_s)\dif s>\lambda\phi(1/(\lambda R))\delta\right\}\\
&\quad\subset\left\{\mathrm{e}^{M^\lambda_{t\wedge\tau_R}}\geq \mathrm{e}^{-\lambda\eps+\lambda\phi(1/(\lambda R))\delta}\right\},
\end{align*}
which then implies the result by Chebyschev's inequality and letting $\lambda=\frac{1}{\eps}$.
\end{proof}

\subsection{A Norris' type lemma}
Let $N(t,\dif y)$ be the Poisson random measure associated with $L_t=W_{S_t}$, i.e.,
$$
N(t,\Gamma)=\sum_{s\leq t}1_\Gamma(L_s-L_{s-}),\ \ \Gamma\in\sB(\mR^d).
$$
Let $\tilde N(t,\dif y)$ be the compensated Poisson random measure of $N(t,\dif y)$, i.e.,
$$
\tilde N(t,\dif y)=N(t,\dif y)-t\nu_L(\dif y),
$$
where $\nu_L$ is the L\'evy measure of $L_t$ given by (\ref{EW2}). By L\'evy-It\^o's decomposition, we have
\begin{align}
L_t=W_{S_t}=W_{\vartheta t}+\int_{|u|\leq 1}y\tilde N(t,\dif y)+\int_{|u|>1}y N(t,\dif y),\label{EW6}
\end{align}
where we have used that for any $0<r<R<\infty$,
$$
\int_{r<|y|\leq R}y\nu_L(\dif y)=0.
$$ 
Notice that $(W_{\vartheta t})_{t\geq 0}$, $\(\int_{|u|\leq 1}y\tilde N(t,\dif y)\)_{t\geq 0}$ and $\(\int_{|u|>1}y N(t,\dif y)\)_{t\geq 0}$ are independent.

Recall the following result about the exponential estimate of martingales (cf. \cite[p.352, (A.5)]{Nu} and \cite[Lemma 1]{Ca}).
\bl\label{Lemma1} Let $\delta, R,\eta, T>0$. 
\begin{enumerate}[(i)]
\item Let $M_t$ be a continuous square integrable martingale, then
$$
P\left\{\sup_{s\in[0,T]}|M_s|\geq\delta; \<M\>_T<\eta\right\}\leq 2\exp\left\{-\frac{\delta^2}{2\eta}\right\}.
$$
\item Let $f_t(y)$ be a bounded $\sF_t$-predictable process with bound $R$, then
$$
P\left\{\sup_{t\in[0,T]}\left|\int^t_0\!\!\!\int_{\mR^d}f_s(y)\tilde N(\dif s,\dif y)\right|\geq\delta,
\int^T_0\!\!\!\int_{\mR^d}|f_s(y)|^2\nu_L(\dif y)\dif s<\eta\right\}\leq 2\exp\left(-\frac{\delta^2}{2(R\delta+\eta)}\right).
$$
\end{enumerate}
\el
The following lemma is contained in the proof of Norris' lemma (cf. \cite[p.137]{Nu} and \cite{Zh1}).
\bl\label{Le4}
For $T>0$, let $f$ be a bounded measurable $\mR^d$-valued function on $[0,T]$. Assume that for some $\eps<T$ and $x\in\mR^d$,
\begin{align}
\int^T_0\left|x+\int^t_0f_s\dif s\right|^2\dif t\leq\eps^3.\label{ET4}
\end{align}
Then we have
$$
\sup_{t\in[0,T]}\left|\int^t_0f_s\dif s\right|\leq 2(1+\|f\|_\infty)\eps.
$$
\el

We now prove the following Norris' type lemma (cf. \cite{No, Nu, Ca, Zh1}).
\bl\label{Lemma2}
Let $Y_t=y+\int^t_0\beta_s\dif s$ be an $\mR^d$-valued process, where $\beta_t$ takes the following form:
$$
\beta_t=\beta_0+\int^t_0\gamma_s\dif s+\int^t_0 Q_s\dif W_{\vartheta s}+\int^t_0\!\!\!\int_{\mR^d}g_s(y)\tilde N(\dif s,\dif y),
$$
where $\gamma_t:\mR_+\to\mR^d$, $ Q_t:\mR_+\to\mR^d\times\mR^d$ and $g_t(y):\mR_+\times\mR^d\to\mR^d$ are three left continuous $\sF_t$-adapted processes.
Suppose that for some left continuous $\sF_t$-adapted $\mR_+$-valued process $\alpha_t$,
\begin{align}
|g_t(y)|\leq \alpha_t(1\wedge|y|).
\end{align}
Then there exists a constant $C\geq1$ such that for any $t\in(0,1)$, $\delta\in(0,\frac{1}{3})$, $\eps\in(0,t^3)$ and $R\geq 1$,
\begin{align}
P\left\{\tau_R>t, \int^t_0|Y_s|^2\dif s<\eps, \int^t_0|\beta_s|^2\dif s\geq 9R^2\eps^{\delta}\right\}\leq
4\exp\left\{-\frac{\eps^{\delta-\frac{1}{3}}}{CR^4}\right\},\label{EW5}
\end{align}
where
$$
\tau_R:=\inf\Big\{t\geq 0: |\beta_t|+|\gamma_t|+| Q_t|+\alpha_t>R\Big\}.
$$
\el
\begin{proof}
Let us define
$$
h_t:=\int^t_0\beta_s\dif s,\quad M^{\mathrm{c}}_t:=\int^t_0\<h_s,  Q_s\dif W_{\vartheta_s}\>,\quad M^{\mathrm{d}}_t:=\int^t_0\!\!\!\int_{\mR^d}\<h_s, g_s(y)\>\tilde N(\dif s,\dif y),
$$
and
\begin{align*}
&E_1:=\left\{\int^t_0|Y_s|^2\dif s<\eps\right\},\ \ E_2:=\left\{\sup_{s\in[0,t]}|h_s|\leq 4R\eps^{\frac{1}{3}}\right\},\\
&E_3:=\left\{\<M^{\mathrm{c}}\>_t\leq C_0R^4\eps^{\frac{2}{3}}\right\},\ \ E_4:=\left\{\sup_{s\in[0,t]}|M^{\mathrm{c}}_s|\leq\frac{\eps^\delta}{2}\right\},\\
&E_5:=\left\{\<M^{\mathrm{d}}\>_t\leq C_1R^4\eps^{\frac{2}{3}}\right\},\ \ E_6:=\left\{\sup_{s\in[0,t]}|M^{\mathrm{d}}_s|\leq\frac{\eps^\delta}{2}\right\},\\
&E_7:=\left\{\int^t_0|\beta_s|^2\dif s<9R^2\eps^\delta\right\},
\end{align*}
where $C_0$ and $C_1$ are two constants determined below.

First of all, by Lemma \ref{Le4}, one sees that for $\eps<T^3$,
\begin{align}
\{\tau_R>t\}\cap E_1\subset\{\tau_R>t\}\cap  E_2\subset \{\tau_R>t\}\cap E_3\cap E_5,\label{EK1}
\end{align}
where the second inclusion is due to
\begin{align*}
\<M^{\mathrm{c}}\>_t=\int^t_0|\<h_s,  Q_s\vartheta\>|^2\dif s\leq 
(4R)^2R^2|\vartheta|^2\eps^{\frac{2}{3}}=:C_0R^4\eps^{\frac{2}{3}}
\end{align*}
and
\begin{align*}
\<M^{\mathrm{d}}\>_t=\int^t_0\!\!\!\int_{\mR^d}|\<h_s, g_s(y)\>|^2\nu_L(\dif y)\dif s\leq 
(4R)^2R^2\left(\int_{\mR^d}1\wedge|y|^2\nu_L(\dif y)\right)\eps^{\frac{2}{3}}=:C_1R^4\eps^{\frac{2}{3}}.
\end{align*}
On the other hand, by integration by parts formula, we have
\begin{align*}
\int^t_0|\beta_s|^2\dif s=\int^t_0\<\beta_s, \dif h_s\>
=\<\beta_t, h_t\>-\int^t_0\<h_s, \gamma_s\>\dif t-M^{\mathrm{c}}_t-M^{\mathrm{d}}_t.
\end{align*}
From this, one sees that on $\{\tau_R>t\}\cap E_2\cap E_4$,
$$
\int^t_0|\beta_s|^2\dif s\leq 4R^2\eps^{\frac{1}{3}}(1+t)+\eps^\delta
\leq(8R^2+1)\eps^\delta\leq 9R^2\eps^\delta.
$$
This means that
$$
\{\tau_R>t\}\cap E_2\cap E_4\cap E_6\subset \{\tau_R>t\}\cap E_7,
$$
which together with (\ref{EK1}) gives
\begin{align*}
&\{\tau_R>t\}\cap E_1\cap E_7^c\subset \{\tau_R>t\}\cap E_1\cap \(E^c_4\cup E^c_6\)\\
&\quad\subset \(\{\tau_R>t\}\cap E_3\cap E^c_4\)\cup\(\{\tau_R>t\}\cap E_5\cap E_2\cap E^c_6\).
\end{align*}
Thus, by Lemma \ref{Lemma1} we have
\begin{align*}
P\(\{\tau_R>t\}\cap E_1\cap E_7^c\)&\leq P\(E_3\cap E^c_4\)+P\(\{\tau_R>t\}\cap E_2\cap E_5\cap E^c_6\)\\
&\leq2\exp\left\{-\frac{\eps^{2(\delta-\frac{1}{3})}}{8C_0R^4}\right\}+2\exp\left(-\frac{\eps^{2\delta}}{8(4R^2\eps^{\frac{1}{3}+\delta}+C_1R^4\eps^{\frac{2}{3}})}\right)\\
&\leq 2\exp\left\{-\frac{\eps^{\delta-\frac{1}{3}}}{8C_0R^4}\right\}+2\exp\left\{-\frac{\eps^{\delta-\frac{1}{3}}}{8(4+C_1)R^4}\right\},
\end{align*}
and (\ref{EW5}) follows by choosing $C:=8(C_0\vee(4+C_1))$.
\end{proof}
\subsection{Malliavin's calculus}
In this subsection we recall some basic notions and facts about the Malliavin calculus (cf. \cite{Ku-St1, Ma, Nu}). Let $\mU$ be a real separable Hilbert space.
Let $\sC(\mU)$ be the class of all $\mU$-valued smooth cylindrical functionals on $\Omega$ with the form:
$$
F=\sum_{i=1}^mf_i(W(h_1),\cdots,W(h_n)) u_i,
$$  
where $f_i\in C^\infty_p(\mR^n)$, $u_i\in\mU$, $h_1,\cdots, h_n\in\mH$ and
$$
W(h)=\int^\infty_0 h_s\dif W_s.
$$
The Malliavin derivative of $F$ is defined by
$$
DF:=\sum_{i=1}^m\sum_{j=1}^n(\p_j f_i)(W(h_1),\cdots,W(h_n))u_i\otimes h_j\in\mU\otimes\mH.
$$
By an iteration argument, for any $k\in\mN$, the higher order Malliavin derivative $D^kF$ of $F$ can be defined as a random variable in $\mU\otimes\mH^{\otimes k}$. 
It is well known that the operator $(D^k,\sC(\mU))$ is closable from $L^p(\Omega;\mU)$ to $L^p(\Omega;\mU\otimes\mH^{\otimes k})$ for each $p\geq 1$ 
(cf. \cite[p.26, Proposition 1.2.1]{Nu}). For every  $p\geq 1$ and $k\in\mN$, we introduce a norm on $\sC(\mU)$ by
$$
\|F\|_{k,p}:=\left(\mE|F|^p+\sum_{l=1}^k\mE\(\|D^lF\|^p_{\mH^{\otimes l}}\)\right)^{\frac{1}{p}}.
$$
The Wiener-Sobolev space $\mD^{k,p}(\mU)$ is defined as the closure of $\sC(\mU)$ with respect to the above norm.
Below we shall simply write
$$
\mD^\infty(\mU):=\cap_{m\in\mN,p\geq 1}\mD^{m,p}(\mU)
$$
and
$$
\mD^{k,p}:=\mD^{k,p}(\mR^d),\ \ \mD^\infty:=\mD^\infty(\mR^d).
$$
The dual operator $D^*$ of $D$ (also called divergence operator) is defined by
$$
\mE \<DF,U\>_\mH=\mE (FD^*U),\ \ U\in \mathrm{Dom}(D^*)=\mD^{1,2}(\mH).
$$
The following Meyer's inequality holds (cf. \cite[p.75, Proposition 1.5.4]{Nu}). For any $p>1$ and $U\in \mD^{1,p}(\mH)$,
\begin{align}
\|D^* U\|_p\leq C_p \|U\|_{1,p}.\label{Me}
\end{align}
Let $F=(F^1,\cdots, F^d)$ be a random vector in $\mD^{1,2}$. The Malliavin covariance matrix of $F$ is defined by
$$
(\Sigma_F)_{ij}:=\<DF^i, DF^j\>_\mH
$$

The following theorem about the criterion that a random vector admits a smooth density in the Malliavin calculus can be found in \cite[p.100-103]{Nu}.
\bt\label{Th2}
Assume that $F=(F^1,\cdots, F^d)\in\mD^\infty$ is a smooth Wiener functional and satisfies that for all $p\geq 2$,
$$
\mE[(\det\Sigma_F)^{-p}]<\infty.
$$
Let $G\in\mD^\infty$ and $\varphi\in C^\infty_p(\mR^d)$. Then for any multi-index $\alpha=(\alpha_1,\cdots,\alpha_m)\in \{1,2,\cdots,d\}^m$, 
$$
\mE[\p_\alpha\varphi(F)G]=\mE[\varphi(F)H_{\alpha}(F,G)],
$$
where $\p_\alpha=\p_{\alpha_1}\cdots\p_{\alpha_m}$, and $H_{\alpha}(F,G)$ are recursively defined by
\begin{align}
H_{(i)}(F,G)&:=\sum_{j}D^*\(G(\Sigma^{-1}_F)_{ij}DF^j\),\label{Me1}\\
H_\alpha(F,G)&:=H_{(\alpha_m)}(F,H_{(\alpha_1,\cdots,\alpha_{m-1})}(F,G)).\no
\end{align}
As a consequence, for any $p\geq 1$, there exists $p_1,p_2,p_3>1$ and $n_1,n_2\in\mN$ such that 
\begin{align}
\|H_\alpha(F,G)\|_p\leq C\|(\det\Sigma_F)^{-1}\|^{n_1}_{p_1}\|DF\|_{m,p_2}^{n_2}\|G\|_{m,p_3}.\label{Esu}
\end{align}
In particular, the law of $F$ possesses an infinitely differentiable density $\rho\in\cS(\mR^d)$.
\et
About the estimate of the density, we recall the following result from Kusuoka-Stroock \cite[Theorem 1.28]{Ku-St1}.
\bt
In the situation of Theorem \ref{Th2}, for any $q>d$, there exists a constant $C=C(q,d)>0$ such that for any $\psi\in C^\infty(\mR^d)$,
\begin{align}
\sup_{y\in\mR^d}|\psi(y)\rho(y)|\leq C\|\psi(F)\|_q^{1-\frac{d}{q}}\left(\sum_i\|H_{(i)}(F,1)\|_q\right)^{d-\frac{d}{q}}
\left(\sum_i\|H_{(i)}(F,\psi(F))\|_q\right)^{\frac{d}{q}},\label{EO1}
\end{align}
provided that the right hand side is finite.
\et

\section{Exponential moment estimate  for SDEs driven by $W_{S_t}$}

In this section, we mainly prove some estimates about the exponential moments for the solutions of
SDEs with non-Lipschitz coefficients. Consider the following SDE driven by $W_{S_t}$:
\begin{align}
\dif X_t=b(X_t)\dif t+A\dif W_{S_t},\ \ X_0=x,\label{EQ99}
\end{align}
where $b:\mR^d\to\mR^d$ is a smooth function and $A=(a_{ij})$ is a constant $d\times d$ matrix.

Recall that a $C^2$-function $H:\mR^d\to\mR^+$ is called a Lyapunov function if
\begin{align}
\lim_{|x|\to\infty}H(x)=\infty.\label{E00}
\end{align}
We assume that for some Lyapunov function $H$ and $\kappa_1,\kappa_2,\kappa_3\geq0$,
\begin{align}
b(x)\cdot\nabla H (x)&\leq \kappa_1 H(x),\label{EE1}
\end{align}
and for all $k=1,\cdots,d$,
\begin{align}
&\Big|\sum_i\p_i H(x)a_{ik}\Big|^2\leq \kappa_2H (x),\label{EE2}\\
&\sum_{ij}\p_i\p_j H(x)a_{ik}a_{jk}\leq \kappa_3.\label{EE3}
\end{align}
Moreover, we also assume the following local Lipschitz condition: for any $R>0$ and all $x,y\in\mR^d$ with $H(x),H(y)\leq R$,
\begin{align}
|b(x)-b(y)|&\leq C_R|x-y|.\label{EE11}
\end{align}

It should be noticed that stochastic differential equation (\ref{EQ99}) can not be solved by using $Y_t=X_t-AW_{S_t}$ to transform (\ref{EQ99}) into an ordinary differential equation with {\it time-dependent} coefficients, since the above conditions are not invariant under this transform. Moreover, a direct application of It\^o's formula seems not work because of
the nonlocal feature of L\'evy processes. 

The main aim of this section is to prove the following estimate.
\bt\label{Th1}
Assume (\ref{EE1})-(\ref{EE11}).  For any initial value $x\in\mR^d$, there exists a unique c\`adl\`ag $\sF_t$-adapted process $t\mapsto X_t$
solving equation (\ref{EQ99}), and for all $t\geq 0$,
\begin{align}
\mE\left[\exp\left\{\frac{2\sup_{s\in[0,t]}H(X_s)}{\e^{\kappa_1 t}(\kappa_2|S_t|+1)}\right\}\right]
\leq C_{\kappa_2,\kappa_3}\e^{H(x)}.\label{EE90}
\end{align}
Moreover, if we also assume {\bf (H$^{2}_{\nu_S}$)}, then for any $p\geq1$ and $t\geq 0$,
\begin{align}
\mE\left[\exp\left\{p\sup_{s\in[0,t]}H(X_s)^{\frac{1}{2}}\right\}\right]
\leq C_{\kappa_1,\kappa_2,\kappa_3,p,t}\e^{H(x)}.\label{EE92}
\end{align}
\et

\begin{proof}
First of all, by (\ref{EE11}) it is standard to prove the existence and uniqueness of local solutions
for equation (\ref{EQ99}). Our main aim is to prove the apriori estimate (\ref{EE90}).
We shall use the approximation argument as used in \cite{Zh1}.

For fixed $\ell\in\mS$, consider the following SDE driven by discontinuous martingale $W_{\ell_t}$:
\begin{align}
\dif X^\ell_t=b(X^\ell_t)\dif t+A\dif W_{\ell_t}, X^\ell_0=x.\label{EE9}
\end{align}
Clearly, it suffices to prove that there exists a unique c\`adl\`ag function $t\mapsto X^\ell_t$
solving equation (\ref{EE9}), and for all $t\geq 0$,
\begin{align}
\mE\left[\exp\left\{\frac{2\sup_{s\in[0,t]}H(X^\ell_s)}{\e^{\kappa_1 t}(\kappa_2|\ell_t|+1)}\right\}\right]
\leq C_{\kappa_2,\kappa_3}\e^{H(x)}.\label{EE93}
\end{align}
Below, for the simplicity of notations, we drop the superscripts ``$\ell$'', and divide the proof into four steps.

{\bf (Step 1).} Let us first consider the case that each component of $\ell$ is absolutely
continuous and strictly increasing.
By It\^o's formula, (\ref{EE1}) and (\ref{EE3}), we have
\begin{align}
\e^{-\kappa_1 t} H(X_t)&= H(x)+\int^t_0\e^{-\kappa_1 s}(b\cdot\nabla H-\kappa_1  H)(X_s)\dif s
+\int^t_0\e^{-\kappa_1 s}\<\nabla H(X_s),A\dif W_{\ell_s}\>\no\\
&\quad+\frac{1}{2}\sum_{ijk}\int^t_0\e^{-\kappa_1 s}\p_i\p_j H(X_s)a_{ik}a_{jk}\dif \ell^k_s\no\\
&\leq H(x)+\int^t_0\e^{-\kappa_1 s}\<\nabla H(X_s),A\dif W_{\ell_s}\>+\frac{\kappa_3}{2}|\ell_t|.\label{EQ6}
\end{align}
For $R>0$, define the stopping time
$$
\tau_R:=\inf\{t\geq 0: |X_t|\geq R\}.
$$
Taking expectations for both sides of (\ref{EQ6}), we obtain that for all $t>0$,
$$
\mE \Big(\e^{-\kappa_1 (t\wedge\tau_R)}H(X_{t\wedge\tau_R})\Big)\leq H(x)+\frac{\kappa_3}{2}|\ell_t|.
$$
This implies by (\ref{E00}) that
\begin{align}
\lim_{R\to\infty}\tau_R=\infty.\label{EW8}
\end{align}

{\bf (Step 2).} Write for $\lambda>0$,
$$
M^\lambda_t:=\lambda\int^t_0\e^{-\kappa_1 s}\<\nabla H(X_s), A\dif W_{\ell_s}\>.
$$
Then by (\ref{EQ6}), we have
\begin{align}
\exp\{\lambda \e^{-\kappa_1 t} H(X_t)\}\leq\exp\left\{\lambda H(x)+\tfrac{\lambda\kappa_3}{2}|\ell_t|\right\}
\exp\left\{M^\lambda_t\right\}.\label{EE4}
\end{align}
Notice that $t\mapsto M^\lambda_t$ is a continuous local martingale with covariance
\begin{align}
\<M^\lambda\>_t:=\lambda^2\sum_k\int^t_0\Big|\e^{-\kappa_1 s}\sum_i\p_i H(X_s)a_{ik}\Big|^2\dif \ell^k_s
\stackrel{(\ref{EE2})}{\leq} \lambda^2\kappa_2G_t|\ell_t|,\label{EE6}
\end{align}
where
\begin{align}
G_t:=\sup_{s\in[0,t]}(\e^{-\kappa_1 s}H(X_s)).\label{EE5}
\end{align}
By Novikov's criterion (cf. \cite{Pr}), one knows that
$$
t\mapsto\exp\Big\{M^\lambda_{t\wedge\tau_R}-\tfrac{1}{2}\<M^\lambda\>_{t\wedge\tau_R}\Big\}
\mbox{ is a continuous exponential martingale},
$$
and by Doob's inequality about positive submartingales and H\"older's inequality,
\begin{align*}
&\mE\exp\left\{\sup_{s\in[0,t]}M^\lambda_{s\wedge\tau_R}\right\}
\leq2\left(\mE\exp\left\{2M^\lambda_{t\wedge\tau_R}\right\}\right)^{\frac{1}{2}}\\
&\quad\leq2\left(\mE\exp\left\{M^{4\lambda}_{t\wedge\tau_R}-\tfrac{1}{2}\<M^{4\lambda}\>_{t\wedge\tau_R}
\right\}\right)^{\frac{1}{4}}\left(\mE\exp\left\{8\<M^{\lambda}\>_{t\wedge\tau_R}\right\}\right)^{\frac{1}{4}}\\
&\quad=2\left(\mE\exp\left\{8\<M^{\lambda}\>_{t\wedge\tau_R}\right\}\right)^{\frac{1}{4}}.
\end{align*}
Recalling (\ref{EE5}) and by (\ref{EE4}) and (\ref{EE6}), we have
\begin{align*}
\mE\exp\left\{\lambda G_{t\wedge\tau_R}\right\}&\leq\exp\left\{\lambda H(x)
+\tfrac{\lambda\kappa_3}{2}|\ell_t|\right\}\mE\exp\left\{\sup_{s\in[0,t]}M^\lambda_{s\wedge\tau_R}\right\}\\
&\leq 2\exp\left\{\lambda H(x)+\tfrac{\lambda\kappa_3}{2}|\ell_t|\right\}
\left(\mE\left\{8\lambda^2\kappa_2G_{t\wedge\tau_R}|\ell_t|\right\}\right)^{\frac{1}{4}}.
\end{align*}
Thus, if one takes $\lambda=\frac{1}{8(\kappa_2|\ell_t|+1)}$, then
$$
\mE\exp\left\{\frac{G_{t\wedge\tau_R}}{8(\kappa_2|\ell_t|+1)}\right\}
\leq 2^{\frac{4}{3}}\exp\left\{\frac{H(x)}{6(\kappa_2|\ell_t|+1)}+\frac{\kappa_3}{12\kappa_2}\right\}
\leq C_{\kappa_2,\kappa_3}\e^{H(x)}.
$$
Finally, by Fatou's lemma and (\ref{EW8}), letting $R\to\infty$, we get
\begin{align}
\mE\exp\left\{\frac{\sup_{s\in[0,t]}H(X_s)}{8\e^{\kappa_1 t}(\kappa_2|\ell_t|+1)}\right\}
\leq\mE\exp\left\{\frac{\sup_{s\in[0,t]}(\e^{-\kappa_1 s}H(X_s))}{8(\kappa_2|\ell_t|+1)}\right\}
\leq C_{\kappa_2,\kappa_3} \e^{H(x)}.\label{EE7}
\end{align}

{\bf(Step 3).} For general $\ell\in\mS$, let us define the Stelkov's average of $\ell$ by
$$
\ell^n_t:=n\int^{t+1/n}_t\ell_s\dif s+\frac{t}{n}=\int^1_0\ell_{t+s/n}\dif s+\frac{t}{n}.
$$
It is clear that $t\mapsto\ell^n_t$ is absolutely continuous and strictly increasing. Moreover, for each $t>0$,
\begin{align}
\ell^n_t\downarrow \ell_t.\label{G6}
\end{align}
By (\ref{EE7}) one has the following uniform estimate:
\begin{align}
\mE\exp\left\{\frac{\sup_{s\in[0,t]}H(X^{\ell^n}_s)}{8\e^{\kappa_1 t}(\kappa_2|\ell^n_t|+1)}\right\}
\leq C_{\kappa_2,\kappa_3}\e^{H(x)}.\label{EE8}
\end{align}
If we define
$$
\tau^{n}_{R_1}:=\inf\Big\{t\geq 0: H(X^{\ell^n}_t)\geq R_1\Big\}
$$
and
$$
\tau_{R_2}:=\inf\Big\{t\geq 0: H(X^\ell_t)\geq R_2\Big\},
$$
then by (\ref{EE11}) and equation (\ref{EE9}), we have for $t<\tau^{n}_{R_1}\wedge\tau_{R_2}$,
\begin{align*}
|X^{\ell^n}_t-X^{\ell}_t|&\leq\int^t_0|b(X^{\ell^n}_s)-b(X^{\ell}_s)|\dif s
+\|A\|\cdot|W_{\ell^n_t}-W_{\ell_t}|\\
&\leq C_{R_1\vee R_2}\int^t_0|X^{\ell^n}_t-X^{\ell}_t|\dif s+\|A\|\cdot|W_{\ell^n_t}-W_{\ell_t}|,
\end{align*}
which yields by Gronwall's inequality that
\begin{align*}
|X^{\ell^n}_t-X^{\ell}_t|&\leq\|A\|\cdot|W_{\ell^n_t}-W_{\ell_t}|
+\exp\left\{C_{R_1\vee R_2}t\right\}\|A\|\int^t_0|W_{\ell^n_s}-W_{\ell_s}|\dif s.
\end{align*}
Now, for any $\eps>0$, by Chebyschev's inequality and (\ref{EE8}), we have
\begin{align*}
P\left\{|X^{\ell^n}_t-X^{\ell}_t|>\eps,t<\tau_{R_2}\right\}&\leq
P\left\{t\geq\tau^{n}_{R_1}\right\}+P\left\{|X^{\ell^n}_t-X^{\ell}_t|>\eps; t<\tau^{n}_{R_1}\wedge\tau_{R_2}\right\}\\
&\leq P\left\{\sup_{s\in[0,t]}H(X^{\ell^n}_s)\geq R_1\right\}+
P\left\{\|A\|\cdot|W_{\ell^n_t}-W_{\ell_t}|>\frac{\eps}{2}\right\}\\
&\quad+P\left\{\exp\left\{C_{R_1\vee R_2} t\right\}
\|A\|\int^t_0|W_{\ell^n_s}-W_{\ell_s}|\dif s>\frac{\eps}{2}\right\}\\
&\leq \frac{2}{R_1}\sup_n\mE\left(\sup_{s\in[0,t]}H(X^{\ell^n}_s)\right)
+\frac{2\|A\|}{\eps}\mE|W_{\ell^n_t}-W_{\ell_t}|\\
&\quad+\frac{2\exp\left\{C_{R_1\vee R_2}t\right\}\|A\|}{\eps}\int^t_0\mE|W_{\ell^n_s}-W_{\ell_s}|\dif s\\
&\leq \frac{C}{R_1}+\frac{2\|A\|}{\eps}|\ell^n_t-\ell_t|^{\frac{1}{2}}
+\frac{2\exp\left\{C_{R_1\vee R_2} t\right\}\|A\|}{\eps}\int^t_0|\ell^n_s-\ell_s|^{\frac{1}{2}}\dif s,
\end{align*}
which tends to zero by (\ref{G6}) as $n\to\infty$ and $R_1\to\infty$. 

Let $\mQ$ be the set of all rational numbers. By a diagonalization argument, there exists a common subsequence $n_m$ and a null set $N$
such that for all $\omega\notin N$ and $t\in\mQ\cap[0,\tau_{R_2}(\omega)]$,
$$
\lim_{m\to\infty}|X^{\ell^{n_m}}_t(\omega)-X^{\ell}_t(\omega)|=0.
$$
Thus, by Fatou's lemma and (\ref{EE8}), we obtain
\begin{align*}
\mE\exp\left\{\frac{\sup_{s\in[0,t\wedge\tau_{R_2}]}H(X^{\ell}_s)}{8\e^{\kappa_1 t}(\kappa_2|\ell_t|+1)}\right\}
&=\mE\exp\left\{\frac{\sup_{s\in[0,t\wedge\tau_{R_2}]\cap\mQ}H(X^{\ell}_s)}{8\e^{\kappa_1 t}(\kappa_2|\ell_t|+1)}\right\}\\
&=\mE\exp\left\{\sup_{s\in[0,t\wedge\tau_{R_2}]\cap\mQ}\liminf_{m\to\infty}\frac{H(X^{\ell^{n_m}}_s)}
{8\e^{\kappa_1 t}(\kappa_2|\ell^{n_m}_t|+1)}\right\}\\
&\leq\liminf_{m\to\infty}\mE\exp\left\{\frac{\sup_{s\in[0,t\wedge\tau_{R_2}]\cap\mQ}
H(X^{\ell^{n_m}}_s)}{8\e^{\kappa_1 t}(\kappa_2|\ell^{n_m}_t|+1)}\right\}\\
&\leq C_{\kappa_2,\kappa_3}\exp\left\{H(x)\right\}.
\end{align*}
Finally, letting $R_2\to\infty$, we obtain (\ref{EE93}).

{\bf(Step 4).} As for (\ref{EE92}), by Young's inequality we have
$$
pH^{\frac{1}{2}}(X_s)\leq\frac{H(X_s)}{16\e^{\kappa_1 t}(\kappa_2|S_t|+1)}+C_p(\kappa_2|S_t|+1).
$$
By H\"older's inequality and (\ref{EE90}), we have
$$
\mE\left[\exp\left\{p\sup_{s\in[0,t]}H(X_s)^{\frac{1}{2}}\right\}\right]
\leq C_{\kappa_2,\kappa_3}\e^{H(x)}\left(\mE \e^{2C_p(\kappa_2|S_t|+1)}\right)^{\frac{1}{2}}
\leq C_{\kappa_1,\kappa_2,\kappa_3,p,t}\e^{H(x)},
$$
where the second inequality is due to {\bf (H$^{2}_{\nu_S}$)}.
\end{proof}

\section{Malliavin Covariance Matrix}

In the sequel,  in addition to (\ref{EE1})-(\ref{EE3}), we also assume that for any $m\in\mN_0$ and some $q_m\geq 0$,
\begin{align}
|\nabla^m b(x)|\leq C(H(x)^{q_m}+1),\label{KK}
\end{align}
where 
$$q_1\in[0,\tfrac{1}{2}].$$
By Theorem \ref{Th1}, it is easy to see that 
\begin{align*}
w&\mapsto X_t(x,w,\ell)\in \mD^\infty(\mR^d),\\ 
x&\mapsto X_t(x,w,\ell)\in C^\infty(\mR^d).
\end{align*}
Let $J_t=J_t(x)=\nabla X_t(x)$ be the derivative matrix of $X_t(x)$ with respect to $x$. Then $J_t$ satisfies
\begin{align}
J_t=I+\int^t_0 \nabla b(X_s)\cdot J_s\dif s.\label{Ep2}
\end{align}
Let $K_t$ be the inverse matrix of $J_t$. Then $K_t$ satisfies
\begin{align}
K_t=I-\int^t_0K_s\cdot \nabla b(X_s)\dif s.\label{Ep3}
\end{align}

We prepare the following basic estimates for later use.
\bl\label{Le11}
Assume {\bf (H$^{2}_{\nu_S}$)}. For $x\in\mR^d$, let $X_t(x)$ be the solution of SDE (\ref{EQ99}).
\begin{enumerate}[(i)]
\item
For any $p\geq 1$, there exists a constant $C_p>0$ such that for all $t\in[0,1]$,
\begin{align}
\mE \left(\sup_{s\in[0,t]}|S_s|^p\right)\leq C_p t.\label{EUY1}
\end{align}
\item There exists a constant $C_x>0$ such that for all $t\in[0,1]$ and $\eps>0$,
\begin{align}
P\left\{\sup_{s\in[0,t]}|X_s(x)-x|>\eps\right\}\leq\frac{C_x t}{\eps^2}.\label{EUY2}
\end{align}
\item For any $p\geq2$, there exists a constant $C_{p,x}>0$ such that for all $t\in[0,1]$,
\begin{align}
\mE\left(\sup_{s\in[0,t]}|J_s(x)-I|^p\right)+\mE\left(\sup_{s\in[0,t]}|K_s(x)-I|^p\right)\leq C_{p,x} t^p.\label{EU1}
\end{align}
\item For any $p\geq2$ and $m,k\in\mN_0$ with $m+k\geq 1$, there exists a constant $C_{p,m,k,x}>0$ such that for all $t\in[0,1]$,
\begin{align}
\mE\left(\sup_{s\in[0,t]}\|D^m\nabla^k X_s(x)\|_{\mH^{\otimes^m}}^p\right)\leq C_{p,m,k,x}
\left\{\begin{aligned}
&1,\quad &m=0,k=1;\\
&t,\quad &m=1,k=0;\\
&t^p,\quad &m+k\geq 2.
\end{aligned}
\right.\label{ET6}
\end{align}
\end{enumerate}
\el
\begin{proof}
(i) By (\ref{EUY}) and {\bf (H$^{2}_{\nu_S}$)}, we can write
$$
S_t=t\left(\vartheta+\int_{\mR^d_+}u\nu_S(\dif u)\right)+\int_{\mR^d_+}u \tilde\mu(t,\dif u)=:t\vartheta'+\int_{\mR^d_+}u \tilde\mu(t,\dif u).
$$
By It\^o's formula, we have
\begin{align*}
|S_t|^p&=p\int^t_0|S_s|^{p-2}\<S_s,\vartheta'\>\dif s+\int^t_0\!\!\int_{\mR^d_+}\(|S_{s-}+u|^p-|S_{s-}|^p\)\tilde\mu(\dif s,\dif u)\\
&\quad+\int^t_0\!\!\int_{\mR^d_+}\(|S_{s-}+u|^p-|S_{s-}|^p-p\<u,S_{s-}\>|S_{s-}|^{p-2}\)\nu_S(\dif u)\dif s.
\end{align*}
Taking expectations and by Young's inequality, we obtain
\begin{align*}
\mE|S_t|^p&\leq p|\vartheta'|\int^t_0\mE|S_s|^{p-1}\dif s+p\int^t_0\!\!\int_{\mR^d_+}|u|\((|S_{s-}|+|u|)^{p-1}+|S_{s-}|^{p-1}\)\nu_S(\dif u)\dif s\\
&\leq C_p\int^t_0\mE|S_s|^{p-1}\dif s+C_p t\int_{\mR^d_+}|u|^p\nu_S(\dif u)\leq C_p\int^t_0\mE|S_s|^p\dif s+C_p t,
\end{align*}
which then gives the estimate (\ref{EUY1}) by Gronwall's inequality and that each component of $S_t$ is increasing.

(ii) Noticing that
$$
\sup_{t\in[0,t]}|X_s(x)-x|\leq\int^t_0 |b(X_s(x))|\dif s+\sup_{s\in[0,t]}|W_{S_s}|,
$$
by Chebyschev's inequality, we have
\begin{align*}
P\left\{\sup_{s\in[0,t]}|X_s(x)-x|>\eps\right\}&\leq\frac{2}{\eps^2}\left(t\int^t_0 \mE|b(X_s(x))|^2\dif s+\mE\(\sup_{s\in[0,t]}|W_{S_s}|^2\)\right)\\
&\leq\frac{C}{\eps^2}\left(t\int^t_0 \(\mE|H(X_s(x))|^{q_0}+1\)\dif s+\mE|S_t|\right),
\end{align*}
which yields (\ref{EUY2}) by (\ref{EE92}) and (\ref{EUY1}).

(iii) By (\ref{Ep2}), we have
$$
|J_t-I|\leq\int^t_0 |\nabla b(X_s)|\cdot |J_s-I|\dif s+\int^t_0|\nabla b(X_s)|\dif s.
$$
which yields by Gronwall's inequality that
$$
|J_t-I|\leq \left(\exp\left\{\int^t_0 |\nabla b(X_s)|\dif s\right\} +1\right)\int^t_0|\nabla b(X_s)|\dif s
$$
By (\ref{KK}) with $q_1\in[0,\frac{1}{2}]$ and (\ref{EE92}), we obtain (\ref{EU1}).

(iv) Notice that for $h\in\mH$,
\begin{align}
D_hX_t=\int^t_0\nabla b(X_s)D_hX_s\dif s+Ah_{S_t}.\label{G7}
\end{align}
Let $\{h^n,n\in\mN\}$ be an orthonormal basis of $\mH$. Then
$$
\|DX_t\|_\mH=\left(\sum_n|D_{h^n}X_t|^2\right)^{\frac{1}{2}}\leq \int^t_0|\nabla b(X_s)|\cdot\|DX_s\|_\mH\dif s+\left(\sum_n|Ah^n_{S_t}|^2\right)^{1/2}.
$$
By Gronwall's inequality and (\ref{EQ1}) below, we obtain
$$
\|DX_t\|_\mH\leq\|A\|\cdot |S_t|^{1/2}+\exp\left\{\int^t_0|\nabla b(X_s)|\dif s\right\} \int^t_0\|A\|\cdot|S_s|^{1/2}\dif s,
$$
which, by (\ref{EE92}), (\ref{KK}) with $q_1\in[0,\frac{1}{2}]$, H\"older's inequality and (i), then gives (\ref{ET6}) for $m=1$ and $k=0$. For the general $m$ and $k$,
it follows by similar calculations and induction method.
\end{proof}

\br\label{Re1}
From the above proof, it is easy to see that if $b\in C^\infty_b(\mR^d)$, then 
$C_x$ in (\ref{EUY2}), $C_{p,x}$ in (\ref{EU1}) and $C_{p,m,k,x}$ in (\ref{ET6}) can be independent of $x\in\mR^d$.
\er
We need the following simple formula about the change of variables (cf. \cite{Ku}).
\bl
Let $f:\mR_+\to\mR$ be a bounded measurable function, and $h:\mR_+\to\mR$ an absolutely continuous function with
integrable derivative. Given a c\`adl\`ag increasing function $\ell$, we have
\begin{align}
\int^t_0f_s\dif h_{\ell_s}=\int^{\ell_t}_0f_{\ell^{-1}_s}\dot h_s\dif s,\label{EQ11}
\end{align}
where $\ell^{-1}_t:=\inf\{s\geq 0: \ell_s>t\}$.
\el
\begin{proof}
By definition, it is easy to see that
$$
\ell^{-1}_t>a\Rightarrow t\geq \ell_a,
$$
and
$$
\ell^{-1}_t\leq a\Rightarrow t\leq \ell_a.
$$
Thus, for $0\leq a<b\leq t$ we have
$$
(\ell_a,\ell_b)\subset\{s: \ell^{-1}_s\in(a,b]\}\subset [\ell_a,\ell_b].
$$
Hence,
$$
\int^{\ell_t}_01_{(a,b]}(\ell^{-1}_s)\dot h_s\dif s=\int^{\ell_t}_01_{(\ell_a,\ell_b]}(s)\dot h_s\dif s=h_{\ell_b}-h_{\ell_a}=\int^t_01_{(a,b]}(s)\dif h_{\ell_s}.
$$
In particular, (\ref{EQ11}) holds for step functions. For general bounded measurable $f$, it follows by
a monotone class argument.
\end{proof}

\bl
Let $f,g:\mR_+\to\mR^d$ be two bounded measurable functions, and $\{h^n,n\in\mN\}$ an orthonormal basis of $\mH$. We have
\begin{align}
\sum_n\left(\int^t_0 f_s\cdot\dif h^n_{\ell_s}\right)\left(\int^t_0g_s\cdot\dif h^n_{\ell_s}\right)
=\sum_k\int^t_0f^k_sg^k_s\dif \ell^k_s.\label{EQ1}
\end{align}
\el
\begin{proof}
If we define
$$
\hat f^k_s:=1_{[0,\ell^k_t]}(s)f^k_{(\ell^k_\cdot)^{-1}_s},\ \ s\geq 0,
$$
and let $\hat f_s=(\hat f^1_s,\cdots,\hat f^d_s)$, then by formula (\ref{EQ11}), we have
\begin{align*}
\int^t_0 f_s\cdot\dif h^n_{\ell_s}=\int^\infty_0 \hat f_s\cdot\dot h^n_s\dif s.
\end{align*}
Thus, by Parsavel's equality, the left hand side of (\ref{EQ1}) equals to
\begin{align*}
\sum_n\left(\int^\infty_0 \hat f_s\cdot\dot h^n_s\dif s\right)
\left(\int^\infty_0 \hat g_s\cdot\dot h^n_s\dif s\right)
=\int^\infty_0\hat f_s\cdot\hat g_s\dif s
&=\sum_k\int^{\ell^k_t}_0f^k_{(\ell^k_\cdot)^{-1}_s}g^k_{(\ell^k_\cdot)^{-1}_s}\dif s,
\end{align*}
which then gives (\ref{EQ1}) by (\ref{EQ11}) again.
\end{proof}

The following lemma originally appeared in the proof of \cite[Theorem 3.3]{Ku}.
\bl
Let $(\Sigma_t(x))_{ij}:=\<DX_t^i(x),DX_t^j(x)\>_\mH$ be the Malliavin covariance matrix of $X_t(x)$. We have
\begin{align}
\Sigma_t(x)=J_t(x)\left(\sum_k\int^t_0K_s(x) a_{\cdot k}
(K_s(x)a_{\cdot k})^*\dif S^k_s\right)(J_t(x))^*.\label{Ep4}
\end{align}
\el
\begin{proof}
By (\ref{Ep2}), (\ref{G7}) and the variation of constant formula, we have
$$
D_hX_t(x)=\int^t_0J_tK_s(x) A\dif h_{S_s}.
$$
Let $\{h^n,n\in\mN\}$ be an orthonormal basis of $\mH$. Then,
$$
\Sigma_t(x)=\sum_nD_{h^n}X_t(x)\cdot (D_{h^n}X_t(x))^*
=\sum_n\left(\int^t_0J_t(x)K_s(x) A\dif h^n_{S_s}\right)\cdot\left(\int^t_0J_t(x)K_s(x) A\dif h^n_{S_s}\right)^*,
$$
which in turn gives the formula (\ref{Ep4}) by (\ref{EQ1}).
\end{proof}

\section{Proof of Main Theorem}

In this section we consider the following SDE
\begin{align}
X_t=x+\int^t_0b(X_s)\dif s+AW_{S_ t}\stackrel{(\ref{EW6})}{=}x+\int^t_0b(X_s)\dif s+AW_{\vartheta t}+\int_{\mR^d}Ay\tilde N(t,\dif y),\label{Eq99}
\end{align}
where, in addition to {\bf(H$^{1}_{\nu_S}$)} and {\bf(H$^{2}_{\nu_S}$)}, we assume that $b\in C^\infty(\mR^d)$ and
\begin{enumerate}[{\bf (H$^1_b$)}]
\item For some some Lyapunov function $H$ and $\kappa_1,\kappa_2,\kappa_3>0$,
\begin{align}
b(x)\cdot\nabla H (x)\leq \kappa_1 H(x),\ \ 
\Big|\sum_i\p_i H(x)a_{ik}\Big|^2\leq \kappa_2H (x),\ \ \sum_{ij}\p_i\p_j H(x)a_{ik}a_{jk}\leq \kappa_3,\label{G4}
\end{align}
and for any $m\in\mN_0$, there is a $q_m\geq 0$ such that
\begin{align}
|\nabla^m b(x)|\leq C(H(x)^{q_m}+1),\label{EW22}
\end{align}
where $q_1\in[0,\frac{1}{2}]$.
\item For some $\kappa_4,\kappa_5,\kappa_6>0$,
\begin{align}
|(\nabla b(x+Ay)-\nabla b(x))A|&\leq \kappa_4(1\wedge|y|),\label{EW0}\\
|(\nabla b(x+Ay)-\nabla b(x)-Ay\cdot\nabla^2 b(x))A|&\leq \kappa_5|y|^2,\label{EW000}\\
\inf_{x\in\mR^d}\inf_{|u|=1}\Big(|uA|^2+|u\nabla b(x)A|^2\Big)&=\kappa_6>0\label{EW11}.
\end{align}
\end{enumerate}

The following estimate is the key part for proving the smooth density of $X_t(x)$.
\bl\label{Le10}
Let $\theta$ be given by (\ref{Con1}). For any $p\geq 1$ and $x\in\mR^d$, there exists a constant $C=C(p,d,\theta,x)>0$ such that for all $t\in(0,1]$,
$$
\|(\det \Sigma_t(x))^{-1}\|_p\leq C t^{-\frac{24 d}{\theta}}.
$$
Moreover, if for all $m\in\mN$, $q_m=0$ in (\ref{EW22}), then the above constant $C$ can be independent of $x$.
\el
\begin{proof} We divide the proof into four steps.

{\bf (Step 1).}  Set
\begin{align*}
&Y_t:=uK_tA,\quad\beta_t:=uK_t\nabla b(X_t)A,\quad  Q_t:=uK_t\nabla^2 b(X_t)A, \\
&\gamma_t:=uK_t\left(\sum_ib^i\cdot\p_i\nabla bA-(\nabla b)^2A+\frac{1}{2}\sum_{ijk}(\p_i\p_j\nabla b)a_{ik}a_{jk}\vartheta_k\right)(X_t)\\
&\qquad+uK_t\int_{\mR^d}\Big(\nabla b(X_t+Ay)-\nabla b(X_t)-Ay\cdot\nabla^2 b(X_t)\Big)A\nu_L(\dif y),\\
&g_t(y):=uK_t(\nabla b(X_{t-}+Ay)-\nabla b(X_{t-}))A.
\end{align*}
By equations (\ref{Ep3}), (\ref{Eq99}) and It\^o's formula, one sees that 
$$
Y_t=aA+\int^t_0\beta_s\dif s,
$$ 
and
$$
\beta_t=a\nabla b(x)A+\int^t_0\gamma_s\dif s+\int^t_0 Q_s\dif s+\int^t_0\!\!\!\int_{\mR^d}g_s(y)\tilde N(\dif s,\dif y).
$$
By (\ref{EW22})-(\ref{EW00}), it is easy to see that
$$
|g_t(y)|\leq C|K_t|(1\wedge|y|),
$$
and for some $q>0$,
$$
|\beta_t|+|\gamma_t|+|Q_t|\leq C|K_t|(H(X_t)^q+1).
$$

{\bf (Step 2).} For $R\geq 1$, define the stopping times
$$
\tau_R:=\inf\Big\{t\geq 0: |K_t|\geq R, H(X_t)\geq R\Big\},
$$
and
$$
\tau_0:=\inf\Big\{s\geq 0: |K_s-I|\geq \tfrac{1}{2}\Big\}.
$$
For $\eta\in(0,1)$, set
$$
E^\eps_t:=\left\{\int^t_0|uK_sA|^2\dif s<\eps^\eta\right\},
$$
and for $\delta\in(0,\frac{1}{3})$ and $R\geq 1$,
$$
F^{\eps,R}_t:=\left\{\int^t_0|uK_s\nabla b(X_s)A|^2\dif s<9R^2\eps^{\eta\delta}\right\}.
$$
By Lemma \ref{Lemma2}, there is a constant $C_1>0$ such that for all $0<\eps<t^{3/\eta}\leq 1$ and $R\geq 1$,
\begin{align*}
P(E^\eps_t\cap\{t<\tau_R\})&=P\Big(E^\eps_t\cap (F^{\eps,R}_t)^c\cap\{t<\tau_R\}\Big)+P\(E^\eps_t\cap F^{\eps,R}_t\cap\{t<\tau_R\}\)\\
&\leq4\exp\left\{-\frac{\eps^{\eta(\delta-\frac{1}{3})}}{CR^4}\right\}+P\(E^\eps_t\cap F^{\eps,R}_t\cap\{\tau_0\geq\eps^{\delta\eta}\}\)+P\(\tau_0<\eps^{\delta\eta}\).
\end{align*}
On the other hand, by (\ref{EW11}) we have
\begin{align*}
E^\eps_t\cap F^{\eps,R}_t&\subset\left\{\int^t_0(|uK_sA|^2+|uK_s\nabla b(X_s)A|^2)\dif s<\eps^\eta+9R^2\eps^{\delta\eta}\right\}\\
&\subset\left\{\int^t_0\frac{|uK_sA|^2+|uK_s\nabla b(X_s)A|^2}{|uK_s|^2}|uK_s|^2\dif s<10R^2\eps^{\delta\eta}\right\}\\
&\subset\left\{\kappa_6\int^t_0|uK_s|^2\dif s<10R^2\eps^{\delta\eta}\right\}.
\end{align*}
Since for any $|u|=1$ and $s\in[0,\tau_0]$,
$$
|uK_s|\geq 1-|K_s-I|\geq \tfrac{1}{2},
$$
it is easy to see that for any $\eps<(\frac{ \kappa_6 t}{20 R^2})^{\frac{1}{\delta\eta}}$,
$$
E^\eps_t\cap F^{\eps,R}_t\cap\{\tau_0\geq\eps^{\delta\eta}\}=\emptyset.
$$
Hence, for $\eta\in(0,1)$, $\delta\in(0,\frac{1}{3})$, $R\geq 1$ and $\eps<(\frac{ \kappa_6 t}{20 R^2})^{\frac{1}{\delta\eta}}$,
\begin{align}
P\(E^\eps_t\cap\{t<\tau_R\}\)\leq 4\exp\left\{-\frac{\eps^{\eta(\delta-\frac{1}{3})}}{CR^4}\right\}+P\(\tau_0<\eps^{\delta\eta}\).\label{EW00}
\end{align}

{\bf (Step 3).}  Now, by Lemma \ref{Le1} and (\ref{EW00}), we have
\begin{align}
P\left\{\int^t_0\sum_k|uK_sa_{\cdot k}|^2\dif S^k_s\leq\eps\right\}&\leq
P\left\{\int^t_0\sum_k|uK_sa_{\cdot k}|^2\dif S^k_s\leq\eps, \int^t_0|uK_sA|^2\dif s\geq\eps^\eta; t<\tau_R\right\}\no\\
&\quad+P\left\{\int^t_0|uK_sA|^2\dif s<\eps^\eta;t<\tau_R\right\}+P\(\tau_R\leq t\)\no\\
&\leq \exp\left\{1-\frac{\phi(\eps/R)\eps^\eta}{\eps}\right\}+4\exp\left\{-\frac{\eps^{\eta(\delta-\frac{1}{3})}}{CR^4}\right\}\no\\
&\quad+P\(\tau_0<\eps^{\delta\eta}\)+P\(\tau_R\leq t\),\label{EW3}
\end{align}
where
$$
\eta\in(0,1),\ \ \delta\in(0,\tfrac{1}{3}),\ \ R\geq 1,\ \ \eps<\(\frac{ \kappa_6 t}{20 R^2}\)^{\frac{1}{\delta\eta}}.
$$
For any $p>1$, by Chebyschev's inequality and (\ref{EU1}),  we have
\begin{align}
P(\tau_0<\eps^{\delta\eta})=P\left\{\sup_{s\in(0,\eps^{\delta\eta})}|K_s-I|\geq\tfrac{1}{2}\right\}
\leq 2^p\mE\left(\sup_{s\in(0,\eps^{\delta\eta})}|K_s-I|^p\right)\leq C\eps^{p\delta\eta},\label{EW33}
\end{align}
and by (\ref{EE92}) and (\ref{EU1}),
\begin{align}
P\(\tau_R\leq t\)\leq \frac{1}{R^p}\mE\left(\sup_{s\in[0,t]}\Big(|K_s|+H(X_s)\Big)^p\right)\leq \frac{C}{R^p}.\label{EW333}
\end{align}
Let $\theta$ be given by (\ref{Con1}). If we choose
$$
\eta=\frac{\theta}{2},\ \ \delta=\frac{1}{6}, \  \ R=\eps^{-\frac{\theta}{48}},
$$
then by (\ref{EW3}), (\ref{EW33}) and (\ref{EW333}), we get for any $t\in(0,1)$, $p\geq 1$ and $\eps\in\(0,\(\frac{ \kappa_6 t}{20}\)^{\frac{24}{\theta}}\)$,
\begin{align}
P\left\{\int^t_0\sum_k|uK_sa_{\cdot k}|^2\dif S^k_s\leq\eps\right\}\leq C\eps^p.\label{EW4}
\end{align}

{\bf (Step 4).}  Let
$$
\xi_t:=\inf_{|u|=1}\sum_k\int^t_0|uK_s a_{\cdot k}|^2\dif S^k_s.
$$
Since $S_t$ has finite moments of all orders, by (\ref{ET6}) and a compact argument
(see \cite[p. 133, Lemma 2.3.1]{Nu} for more details), for any $p\geq 1$, there exists a constant $C_0=C_0(p,d,\kappa_6)>0$ such that 
for all $\eps\in(0,C_0t^{\frac{24}{\theta}})$,
$$
P\Big\{\xi_t\leq\eps\Big\}\leq C_p\eps^p.
$$
Hence,  for all $t\in(0,1]$,
\begin{align*}
\mE\(\xi^{-p}_t\)=p\int^\infty_0\lambda^{p-1}P(\xi^{-1}_t\geq\lambda)\dif\lambda
\leq p\int^{C^{-1}_0t^{-\frac{24}{\theta}}}_0\lambda^{p-1}\dif\lambda+C_p\int^\infty_{C^{-1}_0t^{-\frac{24}{\theta}}}\lambda^{-2}\dif\lambda\leq C_pt^{-\frac{24p}{\theta}}.
\end{align*}
The desired estimate now follows by (\ref{EU1}) and noticing that the smallest eigenvalue of a real symmetric matrix $M$ is less than $(\det M)^{\frac{1}{d}}$.
\end{proof}

We are now in a position to prove the following main result of this paper.
\bt\label{Th0}
Under  {\bf(H$^{1}_{\nu_S}$)}, {\bf(H$^{2}_{\nu_S}$)},  {\bf(H$^{1}_b$)} and {\bf(H$^{2}_b$)}, 
the solution $X_t(x)$ of SDE (\ref{Eq99}) admits a smooth density $\rho(t,x,y)$ as a function on $(0,\infty)\times\mR^d\times\mR^d$.
Moreover, we have the following conclusions:
\begin{enumerate}[(i)]
\item For each $t>0$ and $x\in\mR^d$, $\rho(t,x,\cdot)\in\cS(\mR^d)$ and solves the following Fokker-Planck equation:
\begin{align}
\p_t\rho=\cL\rho(t,x,\cdot)+\div(b\rho(t,x,\cdot)),\label{ET8}
\end{align}
where
$$
\cL f(y)=\mathrm{P.V.}\int_{\mR^d}(f(y+Az)-f(y))\nu_L(\dif z)+\frac{1}{2}\sum_{ijk}(\p_i\p_j f)(y)a_{ik}a_{jk}\vartheta_k.
$$
\item There exist constants $\beta_1,\beta_2,\beta_3>0$ only depending on $d$ and $\theta$ such that for all $(t,x,y)\in(0,1]\times\mR^d\times\mR^d$,
\begin{align}
\rho(t,x,y)\leq C_x \left(t^{-\beta_1}\left(1\wedge\frac{t^{\beta_2}}{|x-y|^{\beta_3}}\right)\right),\label{EE0}
\end{align}
where $C_x$ continuously depends on $x$.
\item Suppose that $b\in C^\infty_b(\mR^d)$, then $C_x$ in (\ref{EE0}) can be independent of $x$.
\end{enumerate}
\et
\begin{proof} For $k,m\in\mN$, by the chain rule, we have
\begin{align*}
\nabla^k\mE(\nabla^mf)(X_t(x))
=\sum_{j=1}^k\mE\Big((\nabla^{m+j}f)(X_t(x))G_j(\nabla X_t(x),\cdots,\nabla^k X_t(x))\Big),
\end{align*}
where $\{G_j, j=1,\cdots,k\}$ are real polynomial functions. By Theorem \ref{Th2} and Lemmas \ref{Le11}, \ref{Le10}, 
one finds that there exist  $\gamma_{k,m}>0$ and $C>0$ such that for all $t\in(0,1)$,
\begin{align}
|\nabla^k\mE (\nabla^mf)(X_t(x))|\leq C\|f\|_\infty t^{-\gamma_{k,m}}.\label{ET7}
\end{align}
Now, by Sobolev's embedding theorem (see \cite[pp.102-103]{Nu}), for each $(t,x)\in(0,\infty)\times\mR^d$,
there exists a smooth density $\rho(t,x,\cdot)\in\cS(\mR^d)$. Moreover, 
$$
(x,y)\mapsto \rho(t,x,y)\in C^\infty(\mR^d\times\mR^d).
$$
By It\^o's formula, one sees that $\rho$  satisfies equation (\ref{ET8}).
The smoothness of $\rho(t,x,y)$ with respect to the time variable $t$ follows by equation (\ref{ET8}) and the standard bootstrap argument. 

Now we use a trick of Kusuoka and Stroock \cite{Ku-St2} to prove (\ref{EE0}).
Let $\chi:\mR^d\to[0,1]$ be a smooth cutoff function with $\chi(y)=0$ for $|y|<\frac{1}{2}$
and $\chi(y)=1$ for $|y|>\frac{3}{4}$. For $\eps>0$, define
$$
\chi_\eps(y):=\chi(y/\eps),\ \ \chi_0(y)=1.
$$
By (\ref{EO1}), we have for any $q>d$,
\begin{align}
\|\chi_\eps(\cdot-x)\rho(t,x,\cdot)\|_\infty&\leq C\|\chi_\eps(X_t(x)-x)\|_q^{1-\frac{d}{q}}\left(\sum_i\|H_{(i)}(X_t(x),1)\|_q\right)^{(1-\frac{1}{q})d}\no\\
&\quad\times\left(\sum_i\|H_{(i)}(X_t(x),\chi_\eps(X_t(x)-x))\|_q\right)^{\frac{d}{q}}.\label{EY1}
\end{align}
By (\ref{Me1}), (\ref{Me}) and H\"older's inequality, we have
\begin{align*}
\|H_{(i)}(X_t(x),\chi_\eps(X_t(x)-x))\|_q&\leq\||\nabla\chi_\eps(X_t(x)-x)|\cdot\|\Sigma^{-1}_t\|\cdot \|DX_t\|^2_{\mH}\|_q\\
&\quad+\|\chi_\eps(X_t(x)-x)\cdot\|D\Sigma^{-1}_t\|\cdot \|DX_t\|_{\mH}\|_q\\
&\quad+\|\chi_\eps(X_t(x)-x))\cdot\|\Sigma^{-1}_t\|\cdot |D^*DX_t|\|_q\\
&\leq\||\nabla\chi_\eps(X_t(x)-x)\|_{q_1}\|\Sigma^{-1}_t\|_{q_2}\|DX_t\|^2_{2q_3}\\
&\quad+\|\chi_\eps(X_t(x)-x)\|_{q_1}\|D\Sigma^{-1}_t\|_{q_2}\|DX_t\|_{q_3}\\
&\quad+\|\chi_\eps(X_t(x)-x)\|_{q_1}\|\Sigma^{-1}_t\|_{q_2}\|DX_t\|_{1,q_3},
\end{align*}
where $\frac{1}{q}=\frac{1}{q_1}+\frac{1}{q_2}+\frac{1}{q_3}$. Let $\mathrm{adj}(\Sigma_t)$ be the adjugate matrix of $\Sigma_t$. Observing that
$$
\Sigma^{-1}_t=\det(\Sigma_t)^{-1}\mathrm{adj}(\Sigma_t), \ \ D\Sigma^{-1}_t=\Sigma^{-1}_tD\Sigma_t\Sigma^{-1}_t,
$$
by the definition of $\mathrm{adj}(\Sigma_t)$, we have
$$
\|\Sigma^{-1}_t\|_{q_2}\leq\|\det(\Sigma_t)^{-1}\|_{2q_2}\|\mathrm{adj}(\Sigma_t)\|_{2q_2}\leq C\|\det(\Sigma_t)^{-1}\|_{2q_2}\|DX_t\|^{2(d-1)}_{4(d-1)q_2},
$$
and
\begin{align*}
\|D\Sigma^{-1}_t\|_{q_2}&\leq\|\det(\Sigma_t)^{-1}\|^2_{8q_2}\|\mathrm{adj}(\Sigma_t)\|^2_{8q_2}\|D X_t\|^2_{1, 4q_2}\\
&\leq C\|\det(\Sigma_t)^{-1}\|^2_{8q_2}\|DX_t\|^{4(d-1)}_{16(d-1)q_2}\|D X_t\|^2_{1, 4q_2}.
\end{align*}
On the other hand, by (\ref{EUY2}), we have
$$
\|\chi_\eps(X_t(x)-x)\|_{q_1}\leq P\left\{|X_t(x)-x|\geq\tfrac{\eps}{2}\right\}^{1/q_1}\leq C\frac{t^{1/q_1}}{\eps^{2/q_1}},
$$
and also
$$
\|\nabla\chi_\eps(X_t(x)-x)\|_{q_1}\leq \frac{\|\nabla\chi\|_\infty}{\eps} P\left\{|X_t(x)-x|\geq\tfrac{\eps}{2}\right\}^{1/q_1}\leq C\frac{t^{1/q_1}}{\eps^{1+2/q_1}}.
$$
Combining the above calculations and by Lemma \ref{Le10} and (\ref{ET6}), we obtain
\begin{align*}
\|H_{(i)}(X_t(x),\chi_\eps(X_t(x)-x))\|_q&\leq C\frac{t^{1/q_1}}{\eps^{1+2/q_1}}\cdot t^{-\frac{24d}{\theta}}\cdot t^{\frac{1}{2q_2}+\frac{1}{q_3}}
+C\frac{t^{1/q_1}}{\eps^{2/q_1}}\( t^{-\frac{48d}{\theta}}\cdot t^{\frac{3}{4q_2}+\frac{1}{q_3}}+t^{-\frac{24d}{\theta}}\cdot t^{\frac{1}{2q_2}+\frac{1}{q_3}}\).
\end{align*}
Similarly,  we also have
$$
\|H_{(i)}(X_t(x),1)\|_q\leq C\( t^{-\frac{48d}{\theta}}\cdot t^{\frac{3}{4q_2}+\frac{1}{q_3}}+t^{-\frac{24d}{\theta}}\cdot t^{\frac{1}{2q_2}+\frac{1}{q_3}}\).
$$
In (\ref{EY1}), taking $\eps=0$ and $\eps=|x-y|$ separately,  by careful choices of parameters, we obtain (ii).
As for (iii), it follows by Remark \ref{Re1}.
\end{proof}

Now we can see that Theorem \ref{Main} is an easy application of Theorem \ref{Th0}.

{\it Proof of Theorem \ref{Main}:}  In the situation of Theorem \ref{Main}, we set
$$
A:=\left(
\begin{aligned}
0,&\ 0\\
0,&\ I
\end{aligned}
\right),\quad
b=
\left(
\begin{aligned}
b_1\\
b_2
\end{aligned}
\right).
$$
By (\ref{G1}) and (\ref{G2}), it is easy to see that (\ref{G4}) and (\ref{EW22}) hold. By (\ref{G3}) and (\ref{G5}), one can see that
(\ref{EW0}), (\ref{EW000}) and (\ref{EW11}) hold.

\vspace{5mm}

{\bf Acknowledgements:}

The author is very grateful to Professors Hua Chen, Zhen-Qing Chen and Feng-Yu Wang for their quite useful conversations.
This work is supported by NSFs of China (No. 11271294) and Program for New Century Excellent Talents in University (NCET-10-0654).

\end{document}